\documentclass[leqno,twoside,dvips]{article}
\usepackage[a5paper,layoutwidth=148.5mm,layoutheight=7.77817in,outer=10mm,inner=38pt,tmargin=10mm,bmargin=10mm,footskip=15pt]{geometry}
\usepackage[utf8]{inputenc}
\usepackage[T1]{fontenc}
\usepackage{lmodern}
\usepackage{fix-cm}
\usepackage{textcomp}
\usepackage{textgreek}
\usepackage{eufrak}
\usepackage[mathcal]{euscript}

\usepackage[np,autolanguage]{numprint}
\usepackage{xspace}

\usepackage[hyphens]{url}%
\RequirePackage{etoolbox}
\appto\UrlSpecials{%
  \do\*{\penalty0 }%
}
\providecommand\url[1]{\texttt{\detokenize{#1}}}
\usepackage{natbib}
\let\bibfont=\small
\makeatletter\def\@biblabel#1{#1.}%
\makeatother
\bibpunct{(}{)}{,}{a}{}{,}
\newcommand\urlprefix{}
\newcommand\bbletal{et al.}

\usepackage{amsfonts}
\newcommand\bbland{and}
\newcommand\bbled{ed.}
\newcommand\bbleds{eds.}

\newcommand\bblin{in}
\newcommand\bblchap{ch.}

\newcommand\natexlab[1]{\textit{#1}}
\providecommand{\citenamefont}[1]{#1}
\newcommand{\mockalph}[1]{}
\providecommand{\url}[1]{\texttt{#1}}
\providecommand{\selectlanguage}[1]{\relax}

          \def\bblvol{vol.}

\usepackage[english,dutch,german,british,french]{babel}\frenchsetup{og=«,fg=»}
\newcommand\german[1]{\foreignlanguage{german}{#1}}
\usepackage{lmodern}
\usepackage{fix-cm}

\usepackage[T1]{fontenc}
\usepackage[utf8]{inputenc}
\usepackage{amsmath}
\usepackage{amssymb}
\usepackage{pifont}

\makeatletter
\DeclareFontEncoding{LS1}{}{}
\DeclareFontEncoding{LS2}{}{\noaccents@}
\DeclareSymbolFont{stix-letters}       {LS1}{stix}     {m}{it}
\DeclareSymbolFont{stix-arrows1}       {LS1}{stixsf}   {m} {n}
\DeclareSymbolFont{stix-operators}     {LS1}{stix}     {m} {n}
\DeclareSymbolFont{stix-largesymbols}  {LS2}{stixex}   {m} {n}
\DeclareSymbolFont{stix-bold-operators}{LS1}{stix}     {b} {n}
\DeclareFontSubstitution{LS1}{stix}{m}{n}
\DeclareFontSubstitution{LS2}{stix}{m}{n}
\SetSymbolFont{stix-letters}     {bold}{LS1}{stix}     {b}{it}
\SetSymbolFont{stix-arrows1}     {bold}{LS1}{stixsf}   {b} {n}
\SetSymbolFont{stix-operators}   {bold}{LS1}{stix}     {b} {n}
\SetSymbolFont{stix-largesymbols}{bold}{LS2}{stixex}   {b} {n}

\def\stix@undefine#1{%
    \if\relax\noexpand#1\let#1=\@undefined\fi}
\def\stix@MathSymbol#1#2#3#4{%
    \stix@undefine#1%
    \DeclareMathSymbol{#1}{#2}{#3}{#4}}
\stix@MathSymbol{\stixwedge}{\mathbin}  {stix-operators}{"E1} \let\land=\wedge
\stix@MathSymbol{\stixvee}  {\mathbin}  {stix-operators}{"E2} \let\lor=\vee
\def\bigwedgegras{\DOTSI\bigwedgeop\slimits@}
\def\bigveegras{\DOTSI\bigveeop\slimits@}
\stix@MathSymbol{\stixbigwedgeop}{\mathop}{stix-largesymbols}{"B4}
\stix@MathSymbol{\stixbigveeop}  {\mathop}{stix-largesymbols}{"B5}
\stix@MathSymbol{\stixrightarrow}               {\mathrel}{stix-arrows1}{"99}
\stix@MathSymbol{\stixrightleftarrows}          {\mathrel}{stix-arrows1}{"CB}
\makeatother

\renewcommand\land{\mathbin{\boldsymbol{\stixwedge}}}
\renewcommand\lor{\mathbin{\boldsymbol{\stixvee}}}

\usepackage{graphicx}
\newcommand\imp{\mathrel{\lower.22ex\hbox{\ding{222}}}}
\newcommand\equ{\mathrel{\ooalign{\raise.08ex\hbox{\ding{222}}\cr\cr\hidewidth\lower.52ex\hbox{\reflectbox{\textrm{\ding{222}}}}\hidewidth}}}
\newcommand\tequ{\mathrel{\ooalign{\lower.22ex\hbox{\ding{222}}\cr\cr\hidewidth\hbox{$\Rightarrow$}\hidewidth}}}

\newcommand\IMP{\mathrel{\preccurlyeq}}

\newcommand\rel{\mathrel{\vdash}}

\newcommand\bigland{\mathop{\textstyle\stixbigwedgeop}\nolimits}

\newcommand\lt{\leqslant}

\newcommand*\psc[1]{%
  \setbox0\hbox{$\mathaccent"0362{#1}^H$}%
  \setbox2\hbox{$\mathaccent"0362{\kern0pt#1}^H$}%
  \ifdim\ht0=\ht2 \overline{#1}\else \bar#1\fi
}
\let\pr\psc
\newcommand{\nats}{\mathbb{N}}
\newcommand\down{\mathord{\downarrow}}

\usepackage{amsthm}
\usepackage{thmtools}
\declaretheorem{lemma}
\declaretheorem[sibling=lemma]{theorem}
\declaretheorem[sibling=lemma,style=definition]{definition}

\newcommand{\lbr}{\lbrack\!\lbrack}
\newcommand{\rbr}{\rbrack\!\rbrack}
\newcommand{\sem}[1] {\lbr #1 \rbr}  

\usepackage{ragged2e}
\setbox0\hbox{\ttfamily\ }\newskip\ttspace\ttspace=\wd0 plus\wd0 minus.5\wd0
\newcommand\deutschmark{\footnotemark}
\newcommand\deutschtext[1]{\footnotetext{\selectlanguage{german}#1}}
\newcommand\deutsch[1]{\deutschmark\deutschtext{#1}}

\usepackage{paralist}
\usepackage{titlesec}

\PassOptionsToPackage{urlcolor=magenta,citecolor=blue,anchorcolor=blue,linkcolor=blue,menucolor=cyan}{hyperref}%
\usepackage[pdftitle=,pdfauthor=,pdfdisplaydoctitle=true,pdflang=en-GB,bookmarksopen=false,breaklinks=true,unicode=true,hyperindex=true,pdfstartview=FitH,colorlinks=true,bookmarksnumbered=true]{hyperref}

\ifx\texorpdfstring\undefined\newcommand\texorpdfstring[2]{#1}\fi
\ifx\pdfbookmark[3]{}\fi
\usepackage{breakurl}
\usepackage{doi}
\usepackage{authblk}

\begin{document}
\selectlanguage{english}
\title{An introduction to Lorenzen's `Algebraic and logistic investigations on free lattices' (1951)}
\author[$\dag$]{Thierry Coquand}
\author[$\ddag$]{Henri Lombardi}
\author[$\ddag$]{Stefan Neuwirth}
\affil[$\dag$]{Computer science and engineering department, University of Gothenburg, Sweden}
\affil[$\ddag$]{Université de Franche-Comté, CNRS, UMR 6623, LmB, 25000 Besançon, France\\\url{stefan.neuwirth@univ-fcomte.fr}}

\date{}
\maketitle


\titleformat{\section}[runin]{\normalfont\normalsize\bfseries}{\thesection.\ }{0pt}{}
\titlespacing{\section}{0pt}{\medskipamount}{\fontdimen2\font plus \fontdimen3\font minus \fontdimen4\font}
\titleformat{\subsection}[runin]{\normalfont\normalsize\itshape}{\thesubsection.\ }{0pt}{}
\titlespacing{\subsection}{0pt}{\medskipamount}{\fontdimen2\font plus \fontdimen3\font minus \fontdimen4\font}
\titleformat{\subsubsection}[runin]{\normalfont\normalsize\upshape}{\thesubsubsection.\ }{0pt}{}
\titlespacing{\subsubsection}{0pt}{\medskipamount}{\fontdimen2\font plus \fontdimen3\font minus \fontdimen4\font}

\section{A foreword on Lorenzen and Gödel.}

Lorenzen inscribes his `Algebraic and logistic investigations on free lattices' into Hilbert's programme as updated after Gödel's impossibility theorems: they merely acknowledge a limitation inherent to the axiomatic method and Lorenzen considers them as a hint and not as an obstacle for a firm foundations of mathematics. `\emph{Metamathematik} (\citealp{lorenzen62}) is written especially to show the irrelevance of Gödel's theorem for a methodical set-up of logic' (\citealp[p.~55]{lorenzen87}).

In 1947, Lorenzen becomes aware of Gödel's proof of consistency of the continuum hypothesis, but decides six months later to stop abruptly considering this problem: `The use of the word ``continuum problem'' I have stopped today' (\citealp[see][]{neuwirth21}).

In 1955, Gödel writes a report on Lorenzen on the occasion of his invitation to the Institute for Advanced Study on behalf of Hermann Weyl shortly before his death. In this report, his `Investigations' are considered to be his greatest achievement in logic, whereas his operative mathematics are judged severely; this is coherent with the way Gödel welcomes Lorenzen at his arrival in 1957, as told by Edward \citet[p.~58]{regis87}. However, he writes another report in 1958 at the conclusion of Lorenzen's stay which is quite positive and insightful (\citealp[see][§~4]{coquandlombardineuwirth25}).

\section{Introduction.}\label{sec:introduction}

Lorenzen's `Algebraische und logistische Untersuchungen über freie Verbände' appeared in 1951 in \emph{The Journal of Symbolic Logic}. These `Investigations' have immediately been recognised as a landmark in the history of infinitary proof theory. Their approach and method of proof have not been incorporated into the corpus of proof theory.\footnote{More precisely, Lorenzen proves the admissibility of cut by double induction, on the cut formula and on the complexity of the derivations, without using any ordinal assignment, contrary to the presentation of cut elimination in most standard texts on proof theory.} We propose a translation and this introduction with the intent of giving a new impetus to their reception.

The `Investigations' are best known for providing a constructive proof of consistency for ramified type theory without axiom of reducibility. They do so by showing that it is a part of a trivially consistent `inductive calculus' that describes our knowledge of arithmetic without detour. The proof resorts only to the inductive definition of formulas and theorems.

They propose furthermore a definition of a semilattice, of a distributive lattice, of a pseudocomplemented semilattice, and of a countably complete boolean algebra as deductive calculi, and show how to present them for constructing conservatively the respective free object over a given preordered set. They illustrate that lattice theory is a bridge between algebra and logic for which the construction of an element corresponds to a step in a proof.

We shall describe the history of their reception, which focusses mainly on the \textomega-rule. The fruitfulness of this device is immediately recognised by Kurt Schütte. It triggers the analysis by \citet{ackermann51} of the infinitary inductive definition of the accessibility predicate in \citealt{gentzen36} with the goal of proving transfinite induction up to ordinal terms beyond~\(\varepsilon_0\), which is also taken over by \citet{schuette52}.

This article is an elaboration of the second half of \citealt{coquandneuwirth17}, and it is a sequel to the elaboration \citealt{coquandneuwirth19} of its first half; we have tried to avoid repetitions and to keep this article self-contained; we invite the reader to look up our \citeyear{coquandneuwirth19} paper for more details about the genesis of Lorenzen's work.

\section{The beginnings.}\label{sec:preparation}

Lorenzen will later recall a talk by Gerhard Gentzen on the consistency of elementary number theory in 1937 or 1938 as a trigger for his discovery that the reformulation of ideal theory in lattice-theoretic terms reveals that his `algebraic works [\dots]\ were concerned with a problem that had \emph{formally} the same structure as the problem of consistency of the classical calculus of logic' (letter to \german{Carl Friedrich Gethmann}, \citealp[see][p.~76]{gethmann91}). This explains the title of Lorenzen's article.

\section{The 1944 manuscript.}\label{sec:1944-manuscript}

The preliminary manuscript `\german{Ein halbordnungstheoretischer Widerspruchsfreiheitsbeweis}', published as \citealt{lorenzen20}, contains already the main ideas, the ω-rule, the connection with free constructions in lattice theory, the admissibility of the cut rule, and applies them to a constructive proof of consistency for elementary number theory.

Note that Lorenzen expresses the concept of an admissible rule only in \citealt{lorenzen48b} (see also his self-review \citealp{lorenzen49}, both have been written after \citealt{lorenzen51}) and names it `erlaubt [permitted]'  and then (in \citealp{lorenzen50}) `eliminierbar [eliminable]'  before adopting today's terminology `zulässig [admissible]' in \citealt{lorenzen55}. Petr Sergeevich Novikov expresses this concept contemporaneously (see \citealp{citkin16}).

\section{The 1945 manuscript.}\label{sec:1945-manuscript}

In a letter dated 11 December 1945 (ETH-Bibliothek, Hochschularchiv, Hs.~975:4111), \foreignlanguage{german}{Heinrich Scholz} submits \german{Lorenzen}'s manuscript `\foreignlanguage{german}{Die Widerspruchsfreiheit der klassischen Logik mit verzweigter Typentheorie}'\footnote{`The freedom from contradiction of classical logic with ramified type theory'. A copy of this manuscript can be found in Niedersächsische Staats- und Universitätsbibliothek Göttingen, Cod.\ Ms.~G.~Köthe M~10.} to \foreignlanguage{german}{Paul Bernays} with a `request for judgement'.

\subsection{Stripping away lattice theory.}
\label{sec:stripp-away-latt}

The manuscript begins as follows.
\begin{quotation}
  \noindent The proof of freedom from contradiction undertaken in the sequel originated as an application of a purely algebraic theorem of existence about `free' complete boolean lattices. In the present work, though, I limit myself exclusively to the logistic application and use no algebraic conceptions whatsoever.\deutsch{`Der im folgenden durchgeführte Widerspruchsfreiheitsbeweis ist als eine Anwendung eines rein algebraischen Existenztheorems über ``freie'' vollständige Boole'sche Verbände entstanden. In dieser Arbeit beschränke ich mich jedoch ausschließlich auf die logistische Anwendung und benutze keinerlei algebraische Begriffsbildungen.'}
\end{quotation}
The choice of stripping away lattice theory may be motivated by targeting a public of logicians. In this way, the affinity with the strategy of \citet[IV, \S\,3]{gentzen35b} becomes more visible: the deductive calculus of ramified type theory with the axioms of comprehension, extensionality and infinity, but without the axiom of reducibility, is compared to an inductive calculus that proceeds `without detour'; with respect to Gentzen's calculus, it features an induction rule (compare rule~[4] on p.~
98) with infinitely many premisses, i.e.\ an \textomega-rule in today's terminology. 

\subsection{Formula inductions and theorem inductions.}
\label{sec:form-induct-theor}

Lorenzen emphasises as follows.
\begin{quotation}
  \noindent This proof uses as auxiliary devices only formula inductions vs.\ theorem inductions, i.e.\ the fact that the concept of formula and the concept of theorem is defined inductively. The harmlessness 
  of these auxiliary means seems to me to be even more transparent than the harmlessness of explicit transfinite inductions.\deutsch{`Dieser Beweis benutzt als Hilfsmittel nur Formelinduktionen bzw.\ Satzinduktionen, d.\ h.\ die Tatsache, daß der Formelbegriff und der Satzbegriff induktiv definiert ist. Die Unbedenklichkeit dieser Hilfsmittel scheint mir noch einleuchtender zu sein, als die Unbedenklichkeit expliziter transfiniter Induktionen.'}
\end{quotation}
These inductions establish that the deductive calculus is a part of the inductive calculus 
in section~
7 of the article:
\begin{itemize}
\item[(1a)] the `logical axiom' $c\IMP c$ is proved by formula induction;
\item[(1b)] the axiom of comprehension follows from the construction of a \textlambda-calculus and a rule of constants;
\item[(1c)] the axiom of extensionality results from a formula induction by the aid of two auxiliary rules proved by theorem induction;
\item[(1d)] the axiom of infinity follows from the properties of the order on numbers;
\item[(2a)] the admissibility of the cut rule is proved by a formula induction on the cut formula: if it is a numerical formula, a double theorem induction on the premisses is used; the only difficulties in the induction step result from the copresence of constants and free and bound variables in rules like~[3d] on p.~
  98; as usual, contraction plays an important rôle.
\end{itemize}

\subsection{Bernays' judgement.}
\label{sec:bernays-judgment}

\foreignlanguage{german}{Bernays} is able to appreciate its content on the spot and replies with detailed comments to Scholz on 24 April 1946 (carbon copy, Hs.~975:4112). On 17 April 1946, \foreignlanguage{german}{Lorenzen} writes directly to Bernays (Hs.~975:\allowbreak2947); he gets an answer on 22 May 1946 with the following appreciation.

\begin{quotation}
  \noindent It seems to me that your argumentation accomplishes in effect the desired and that thereby at the same time also a new, methodically more transparent proof of freedom from contradiction for the number-theoretic formalism, as well as for \foreignlanguage{german}{Gentzen}'s subformula theorem\footnote{In the letter of 24 April 1946, Bernays writes more precisely to Scholz `that one also gets a proof for the main theorem of \foreignlanguage{german}{Gentzen}'s ``Investigations into logical deduction'' out of it, if on the one hand one omits the higher axioms~[(1b,\,c,\,d) on p.~
    95] in the deductive calculus, on the other hand one retains from the rules of the inductive calculus (for determining the concept of theorem) only~[[2], [3a--d] on p.~
    98], while one takes also the formula pairs~[$c\IMP c$] as starting theorems for this calculus.'\deutschmark}\deutschtext{`dass man aus ihm auch einen Beweis für den Hauptsatz von Gentzen's ``Untersuchungen über das logische Schliessen'' erhält, indem man einerseits beim deduktiven Kalkul die höheren Axiome~1.)\,b),\,c),\,d) weglässt, andererseits von den Regeln des induktiven Kalkuls (zur Bestimmung des Satzbegriffes) nur 2)\,a)--d) beibehält, während man als Ausgangssätze auch für diesen Kalkul die Formelpaare~$\mathfrak c\subset\mathfrak c$ nimmt.'}
  is provided.
  
  In the circumstance that all this is included in your result shows at the same time the methodical superiority of your method of proof with respect to a proof (that probably did not come to your knowledge) that F. B.~Fitch [\dots]\ gave in 1938, and that also bears on the comparison of the deductive formalism with a system of formulas which is not delimited in a purely operative way; namely, this delimitation is carried out there according to a definition of truth in which the `tertium non datur' (indeed only with respect to the species of natural numbers) is made use of.\footnote{See \citealt{fitch38} and its review \citealt{bernays39}.} By determining your system of comparison according to the idea of a generalisation of \foreignlanguage{german}{Gentzen}'s thought of `deduction without detour', you gain the possibility of applying the constructive proof-theoretic view also in the case of your `inductive calculus', i.e.\ of such an inference system that does not comply with the recursiveness conditions that the customary formalisms fulfil.\deutsch{`Es scheint mir, dass Ihre Beweisführung in der Tat das Gewünschte leistet und dass damit zugleich auch ein neuer, methodisch durchsichtigerer Wf.-Beweis für den zahlentheoretischen Formalismus wie auch ein solcher für Gentzen's Teilformelsatz geliefert wird.

    In dem Umstande, dass alles dies in Ihrem Ergebnis eingeschlossen ist, zeigt sich zugleich die methodische Überlegenheit Ihres Beweisverfahrens gegenüber einem (Ihnen wohl nicht zur Kenntnis gelangten) Beweis, den F. B.~Fitch 1938 für die Widerspruchsfreiheit der verzweigten Typentheorie gegeben hat (im Journal of symb.\ logic, vol.~3, S.~140-149), und der auch auf dem Vergleich des deduktiven Formalismus mit einem Formelsystem beruht, das auf eine nicht rein operative Art abgegrenzt ist; diese Abgrenzung erfolgt nämlich dort im Sinne einer Wahrheitsdefinition, wobei von dem ``tertium non datur'' (allerdings nur demjenigen in Bezug auf die Gattung der natürlichen Zahlen) Gebrauch gemacht wird. Indem Sie Ihr Vergleichssystem gemäss der Idee einer Verallgemeinerung von Gentzen's Gedanken der ``umweglosen Herleitung'' bestimmen, gewinnen Sie die Möglichkeit, die konstruktive beweistheoretische Betrachtung auch im Falle Ihres ``induktiven Kalküls'' anzuwenden, d.~h.\ eines solchen Folgerungssystems, welches nicht den durch die üblichen Formalismen erfüllten Rekursivitätsbedingungen genügt' (Hs.~975:2948).}
\end{quotation}

\subsection{Independence of the axiom of reducibility.}
\label{sec:indep-axiom-reduc}

Lorenzen learns about Fitch's proof of consistency only by this letter. In his answer (dated 7 June 1946, Hs.~975:2949), he explains the lattice-theoretic background of his proof and encloses a manuscript, `\german{Über das Reduzibilitätsaxiom}',\footnote{`On the axiom of reducibility', Hs.~974:\allowbreak149. Another copy of this manuscript is in the \foreignlanguage{german}{Universitätsarchiv Bonn}.} which is a preliminary version of the last section of the published article, in which the axiom of reducibility is shown to be independent. However, Bernays seems to already have received this manuscript with Scholz's letter of 11 December 1945 (see his letter of 24 April 1946).

More precisely, he proves the consistency of the calculus obtained by adding an axiom that expresses that all infinite sets are countable and then shows that the axiom of reducibility is false in this calculus. In fact, as Lorenzen notes, \citet{fitch39} proves this in his framework. These results answer questions raised by \citet[pp.~xiv, xlii--xxliii]{whiteheadrussell25} in the introduction to the second edition of their \emph{Principia mathematica} after Leon \citet{chwistek25}: without the axiom of reducibility, `Cantor's proof that \(2^n>n\) breaks down unless \(n\)~is finite'.

\section{The 1947 manuscript.}\label{sec:1947-manuscript}

\subsection{Restoring the lattice-theoretic part.}
\label{sec:rest-algebr-part}

By a letter dated 21 February 1947, Lorenzen writes to Bernays:
\begin{quotation}
  \noindent After a revision of my proof of freedom from contradiction according to your precious remarks and after addition of an algebraic part, I would like to allow myself to ask you for your intercession for a publication in the Journal of Symbolic Logic.\deutsch{`nach einer Überarbeitung meines Wf.beweises nach Ihren wertvollen Bemerkungen und nach Hinzufügung eines algebraischen Teiles möchte ich mir erlauben, Sie um Ihre Fürsprache zu bitten für eine Veröffentlichung im Journal of symbolic logic' (Hs.~975:2950).}
\end{quotation}
This new draft is a kind of synthesis of `\german{Ein halbordnungstheoretischer Widerspruchsfreiheitsbeweis}' and `\german{{Die Widerspruchsfreiheit der klassischen Logik mit verzweigter Typentheorie}}', or rather a juxtaposition of two parts: the seams remain apparent. However, the introduction now takes into account the added algebraic part. Its first paragraph 
emphasises that lattice theory is relevant for ideal theory with a reference to the reshaping of Krull's Fundamentalsatz for integral domains in lattice-theoretic terms provided by his habilitation (\citealp{lorenzen50b}, see \citealp{neuwirth21b}).

\subsection{Semilattices.}
\label{sec:semilattices}

In the new algebraic part, the construction of free semilattices and free distributive lattices stems in fact from ideal theory. Theorems~1--4 in section~2 of the article (p.~
84) introduce a semilattice as a `single statement entailment relation' and construct the free semilattice over a preordered set (a \emph{preorder} or quasiorder is a reflexive and transitive relation on the set). This approach may be dated back to \citet[§~2]{skolem21}, who constructs the free lattice over a preordered set in the course of studying the decision problem for lattices. It is parallelled in~\citealt{lorenzen52} by the definition of a system of ideals for an arbitrary preordered set~$M$ on which a monoid~$G$ acts by order-preserving operators~$x$: it is a relation satisfying items~1--4 of theorem~1 and furthermore
\[\text{if $a_1,\dotsc,a_n\rel b$, then $xa_1,\dotsc,xa_n\rel xb$}\]
(compare \citealp[§~1C]{coquandlombardineuwirth19}).

\subsection{Distributive lattices.}
\label{sec:distr-latt}

In the same way, theorems~5--8 provide the description of a distributive lattice as a deductive system that has been called since \cite{scott71} an `entailment relation'. This description strikes Bernays as new to him (letter of 3 April 1947, \german{Paul-Lorenzen-Nachlass, Philosophisches Archiv, Universität Konstanz}, PL~1-1-118). His theorem~7 on p.~
85 corresponds in fact to theorem~1 in \citealt{CeCo2000}, obtained independently. This construction is used in \citealt{lorenzen53} for embedding a preordered group endowed with a system of ideals into a lattice-ordered group containing this system (compare \citealt[§~2C]{coquandlombardineuwirth19}).

\subsection{Lorenzen algebras.}
\label{sec:lorenzen-algebras}

Section~3 of the article deals with (finitary) pseudocomplemented semilattices while his 1944 manuscript deals with countably complete ones. We propose to call them `Lorenzen algebras':\footnote{The theory of Lorenzen algebras continues to develop: one can find an account of it by \citet[pp.~99--101]{graetzer11} and by \citet[ch.~3]{chajdahalaskuehr07}.} see definition~\ref{lorenzen-algebra} in our section~\ref{sec:impr-quant-induct}, where we use them in a crucial way for our explanation of an impredicative system in terms of inductive definitions. 

He constructs the free Lorenzen algebra generated by a preordered set as a cut-free sequent calculus (while his 1944 manuscript deals with the countably complete case).\footnote{The existence of the free Lorenzen algebra over a preordered set seems to be unknown in the literature, which considers only the case where the preorder is trivial; in the latter case, the decision problem has been solved by \cite{tamura74}. \citealt{neuwirth15} proposes a streamlined presentation.}

In section~4, he shows how to apply the construction of the free Lorenzen algebra to a simple intuitionistic logical calculus. He emphasises that the decision problem has a positive answer.

\subsection{Boolean algebras.}
\label{sec:boolean-algebras}

Lorenzen proceeds with describing the free countably complete boolean algebra\footnote{In the second paragraph of the introduction, he addresses complete boolean algebras over a preordered set as studied by \citet{macneille37}. The question about the existence of the free complete boolean algebra is usually attributed to \citet{rieger51} and has led to the works of \citet{gaifman64} and \citet{hales64} that provide a negative answer; it may be seen as an anomaly of the set-theoretic approach to actual infinity, because such infinities cause that free complete boolean algebras would be too big. See also the proof of \citet{solovay66} inspired by forcing.} generated by a preordered set as a cut-free infinitary sequent calculus with \textomega-rules~[3.9] and~[3.10] on p.~
92. The main step in the construction is again to prove that the cut rule (which states on the same page 
that \(a_1\lt c\lor b_1\) and \(a_2\land c\lt b_2\) implies \(a_1\land a_2\lt b_1\lor b_2\)) is admissible. 

He sketches this construction, which goes along the same lines as the construction of the free Lorenzen algebra, with one significant difference: in the latter case, he is able to prove contraction (see lemma~(8) on p.~
88), whereas he has to put it into the definition in the former setting (he provides a counterexample on pp.~
92--93). Compare with his 1944 manuscript, where a contraction rule is present (see the comment in our \citeyear[end of~§~2]{coquandneuwirth19}), and with the calculus defined by Per \citet[\S\,30]{PML} for Borel sets, where the problem of contraction is eluded by `identify[ing] sequents which are equal considered as finite sets'.

Let us comment on two aspects of this construction.
\begin{itemize}
\item Lorenzen works systematically with a preorder~\(\lt\) instead of an order and does not quotient with respect to the equivalence relation $a\equiv b$ defined by $a\lt b$ and $b\lt a$. If he did so, he would need to resort to the axiom of choice for defining meets.
\item The universal property corresponding to freeness is proved by parallelling the construction of the free object with the construction of the sought-after morphism: see items~(i)--(iii) on p.~
  89. \citet{hottbook} indicates a way to avoid the use of choice in a constructive setting precisely by defining inductively the free boolean algebra so that its objects and their equalities are defined simultaneously.
\end{itemize}
Compare our section~\ref{sec:impr-quant-induct} for the relevance of these two aspects.

Then Lorenzen shows how to deduce consistency for the logic of ramified type by an iterated construction of free countably complete boolean algebras, starting from a calculus without free variables, along the hierarchy of types.

\section{Toward publication.}\label{sec:toward-publication}

\subsection{Finitary vs.\ constructive logic.}
\label{sec:finit-vs.-constr}

At the end of his letter of 21 February 1947, Lorenzen asks:
\begin{quotation}
  \noindent I beg once again to ask you for your advice---namely, it is not clear to me whether I rightly call the logic used here `finitary' logic.\deutsch{`Ich bitte noch einmal Sie um Ihren Rat fragen zu dürfen --~es ist mir nämlich nicht klar, ob ich die hier benutzte Logik mit Recht ``finite'' Logik nenne' (Hs.~975:2950).}
\end{quotation}
Bernays provides the following answer in his letter of 3 April 1947:
\begin{quotation}
  \noindent When it comes to the methodical standpoint and to the terminology to be used in relation, then it seems advisable to me to keep with the mode chosen by Mister Gentzen, that one speaks of `finitary' reflections only in the narrower sense, i.e.\ relating to considerations that may be formalised in the framework of recursive number theory (possibly with extension of the domain of functions to arbitrary computable functions), that one uses in contrast the expression `constructive' for the appropriate extension of the standpoint of the intuitive self-evidence; by the way, this is employed also by many an American logician in the corresponding sense.

  Your proof of freedom from contradiction cannot, I deem, be a finitary one in the narrower sense. Of course, this would conflict with the Gödel theorem. Actually, a nonfinitary element of your reflection lies in the induction rule of the inductive calculus, which contains indeed a premiss of a more general form.\deutsch{`Was den methodischen Standpunkt und die in Bezug darauf zu verwendende Terminologie betrifft, so erscheint es mir als empfehlenswert, den von Herrn Gentzen gewählten Modus beizubehalten, dass man von `finiten' Betrachtungen nur im engeren Sinne spricht, d.\ h.\ mit Bezug auf Überlegungen, die sich im Rahmen der rekursiven Zahlentheorie (eventuell mit Erweiterung des Funktionenbereiches auf beliebige berechenbare Funktionen) formalisieren lassen, dass man dagegen für die sachgemässe Erweiterung des Standpunktes der anschaulichen Evidenz den Ausdruck `konstruktiv' verwendet; dieser wird übrigens auch von manchen amerikanischen Logikern im entsprechenden Sinn gebraucht.

    Ihr Wf-Beweis kann, so meine ich, kein finiter in dem genannten engeren Sinne sein. Das würde doch dem Gödelschen Theorem widerstreiten. Tatsächlich liegt, so viel ich sehe, ein nicht-finites Element Ihrer Betrachtung in der Induktionsregel des induktiven Kalkuls, welche ja eine Prämisse von allgemeinerer Form enthält' (PL~1-1-118).}
\end{quotation}
In other words, the \textomega-rule does not fit into a formal system, and this explains why Gödel's theorem does not apply here. But \citet[p.~491]{hilbert31} describes the \textomega-rule as a `finitary deduction rule' and this is probably why Lorenzen qualifies his deductions as `finitary'; note also the emphasis of \citet[§~16.1\,1]{gentzen36} that accessibility is a finitary concept (see our section~\ref{sec:ackermann} below). Lorenzen answers as follows on 4 May 1947.
\begin{quotation}
  \noindent Your proposal to call the means of proof not `finitary' but `constructive' acted on me as a sort of redemption. I was sticking so far to the word finitary only to emphasise that these are hilbertian ideas that I am trying to pursue.\deutsch{`Ihr Vorschlag, die Beweismittel nicht ``finit'', sondern ``konstruktiv'' zu nennen, hat wie eine Art Erlösung auf mich gewirkt. Ich habe bisher an dem Wort finit nur festgehalten, um zu betonen, dass es Hilbertsche Ideen sind, die ich fortzuführen versuche' (Hs.~975:2953).}
\end{quotation}

\subsection{Publication.}
\label{sec:toward-publication-1}

Lorenzen prepares another final draft that is very close to the published version.\footnote{Two pages of this draft may be found in the file~OB~5-3b-5; Cod.\ Ms.~G.~Köthe M~10 contains an excerpt of Part~I.} Bernays sends a first series of comments on 1 September 1947 (PL~1-1-112)
and a second series (on a version including the final section on the axiom of reducibility) on 6 February 1949 (PL~1-1-107); the article is submitted to \emph{The Journal of Symbolic Logic} soon afterwards\footnote{See the 
  letter of 27 April 1949 to Alonzo Church, in which Lorenzen thanks him for acknowledging receipt of the manuscript, writes a few words on its history, and proposes Bernays as a referee (Alonzo Church Papers, box~26 folder~4, Manuscripts division, Department of rare books and special collections, Princeton university library).} and published as \citealt{lorenzen51} with date of reception 17 March 1950. In fact, in 1947, Lorenzen already starts his project of language levels which will lead to his operative logic (see \citealp{neuwirth21} on the circumstances of this switch).

\section{Reception.}

\subsection{Early accounts.}

For early accounts of the manuscripts, see \citealt{lorenzen48short}, \citealt{koethe48}, \citealt[\S\,11]{schmidt50b}.
\medskip

Let us present the reactions to Lorenzen's article by subject.

\subsection{The difference with Fitch's proof.}
\label{sec:hao-wang-predicative}

Hao \citet{wang51} writes the review for \emph{The Journal of Symbolic Logic} and tries to compare Lorenzen's approach with Fitch's; see \citealt{coquand14} for a discussion of this review. \citealt{wang54} (p.~252) provides a more accurate comparison.

\subsection{Induction rules instead of transfinite inductions.}
\label{sec:koth-lect-proofs}

Gottfried Köthe and Lorenzen have worked together on lattice theory before World War~II\@. In Spring 1947, they correspond on foundations of mathematics and physics. Köthe is preparing lectures on proofs of consistency up to Lorenzen's to be given in Fall 1947 at Mainz (see Cod.\ Ms.~G.~Köthe G~3).\footnote{See Köthe's letter of 3 September 1947 (PL 1-1-113) and Cod.\ Ms.~G.~Köthe M~10.} In an answer to a letter by Köthe dated 8 June 1947 (PL~1-1-114), Lorenzen writes on 17 June 1947 about his work:
\begin{quotation}
  \noindent The formalisation of proofs of freedom from contradiction that I am striving for is not at all intent on staying inside a transfinite theory of numbers, but uses instead of `transfinite inductions' induction rules like the `formula induction' and `theorem induction' of my proof of freedom from contradiction---these could also be formalised in a theory of numbers with sufficiently large constructible ordinal numbers, but nothing is gained from this.\deutsch{`Die Formalisierung der Wf\,beweise, die ich anstrebe, geht nun gar nicht darauf aus, innerhalb einer transfiniten Zahlentheorie zu bleiben, sondern benutzt statt der ``transfiniten Induktionen'' Induktionsregeln wie die ``Formelinduktion'' und ``Satzinduktion'' meines Wf\,beweises 
  --~-- diese ließen sich zwar auch in einer Zahlentheorie mit genügend großen konstruierbaren Ordinalzahlen formalisieren, dadurch wird aber nichts gewonnen' (Cod.\ Ms.~G.~Köthe M~10).}
\end{quotation}

\subsection{The logical status of the \textomega-rule.}
\label{sec:logic-stat-text}

In his letter dated 9 June 1947, Ackermann would find appropriate that one `would describe somehow the constructive, contentful thinking' in Lorenzen's rule of induction. Four years later, \citet[pp.~368--369]{ackermann52} reflects upon the logical status of the \textomega-rule in the following terms.
\begin{quotation}
  \noindent Indeed, the expression `derivation rule' does not seem entirely appropriate to us. For on the one hand it is a derivation rule with infinitely many premisses we are dealing with, so that no rigorous formalisation of thought is carried out, as the fact that a formal derivation can be given for each of the infinitely many premisses is the result of considerations in terms of content. [\dots{}] It matters to me here now to show that one can add certain basic formulas that have been obtained according to certain principles [\dots{}] without loosing freedom from contradiction.
\end{quotation}
This analysis is the same as that of \citet{hilbert31}. On p.~489, he writes the following about his proof theory and formalisation twenty years earlier: 
\begin{quotation}
  \noindent Deducing in terms of content is superseded by an exterior acting according to rules, namely the use of the deduction scheme and substitution. [\dots{}].

  To the proper thus formalised mathematics comes an in a way new mathematics, a metamathematics, that is necessary to secure the former, in which---contrary to the purely formal ways of deducing of proper ma\-the\-ma\-tics---deducing in terms of content is applied, but only for verifying that the axioms are free from contradiction.

  The axioms and provable assertions, i.e.\ the formulas that arise in this interplay, are the images of the thoughts that have been making up the habitual proceeding of mathematics up to now.
\end{quotation}
On p.~491, he writes:
\begin{quotation}
  \noindent If it is verified that the formula
  \[\mathfrak A(\mathfrak z)\]
  becomes a correct numerical formula every time \(\mathfrak z\) is a presented numeral, then the formula
  \[(x)\mathfrak A(x)\]
  may be put on as starting formula.
\end{quotation}
Thus the \textomega-rule is considered as a rule of metamathematics applied in terms of content and it is ancillary in studying formal systems.

Ackermann discusses a version of Lorenzen's work, most probably the 1945 manu\-script, in a letter to Lorenzen dated 31 April 1950 (sic, PL~1-1-95). \foreignlanguage{german}{\citealt{ackermann53}} provides a version of the consistency of ramified type theory in the context of his type-free logic (see also \citealp{schuette54b}).

\subsection{The \textomega-rule combined with transfinite induction.}
\label{sec:restr-transf-induct}

In an answer dated 4 November 1948 to a letter by Bernays that informs him about Lorenzen's work, Schütte states that Arnold Schmidt has acquainted him with it and writes: `As means of proof going beyond the narrower finitary standpoint, Mister Lorenzen uses inferences with infinitely many premisses, while I (as Gentzen) draw on beginning cases of the transfinite induction.'\deutsch{`Als Beweismittel, die über den engeren finiten Standpunkt hinausgehen, benutzt Herr Lorenzen Schlüsse mit unendlich viel Prämissen, während ich (wie Gentzen) Anfangsfälle der transfiniten Induktion heranziehe' (Hs.~975:4228, note that Szabo translates `Anfangsfälle der' by the epithet `restricted').}

In a letter to Bernays dated 26 August 1949, Schütte writes: `I believe that my investigations are not superfluous  besides those of Lorenzen because with them the required metamathematical means of proof and the connections with the derivability of the formalised transfinite induction are uncovered.'\deutsch{`Ich glaube, daß meine Untersuchungen neben denen von Lorenzen deshalb nicht überflüssig sind, weil die benötigten metamathematischen Beweismittel und die Zusammenhänge mit der Herleitbarkeit der formalisierten transfiniten Induktion dabei aufgedeckt werden' (Hs.~975:4230).}

Schütte writes to Lorenzen on 1 May 1950 in order to acknowledge the latter's priority in implementing the \textomega-rule into proofs of consistency.\footnote{Both are not aware of the work of \citet{novikoff39,novikoff43} in this respect. See \citealt[§~1.2]{mints91}.}
\begin{quotation}
  \noindent [\dots{}]\ I came to know that you had provided already before a proof of freedom from contradiction for a still more general domain, and had arrived at the following result: the cut-eliminability, that with Gentzen had been carried out only in pure logic, may also be transferred to mathematical formalisms, if instead of the inference of complete induction more general schemes of inference with infinitely many premisses are drawn on by extending the concept of derivation so that it may contain infinitely many formulas. This insight gained by you, that appears to me exceptionally important for fundamental research, I have now taken up.\deutsch{`[\dots]\ erfuhr ich, daß Sie schon vorher einen Widerspruchsfreiheitsbeweis für einen noch allgemeineren Bereich erbracht hatten und dabei zu folgendem Ergebnis gekommen waren: Die Schnitt-Eliminierbarkeit, die bei Gentzen nur in der reinen Logik durchgeführt wurde, läßt sich auch auf mathematische Formalismen übertragen, wenn statt des Schlusses der vollständigen Induktion allgemeinere Schlußschemata mit unendlich vielen Prämissen herangezogen werden, indem der Begriff der Herleitung so erweitert wird, dass er unendlich viele Formeln enthalten darf. Diese von Ihnen gewonnene Erkenntnis, die mir außerordentlich wichtig für die Grundlagenforschung zu sein scheint, habe ich nun aufgegriffen' (PL~1-1-95).}
\end{quotation}

\subsection{Semi-formal systems.}
\label{sec:semiformal-systems}

In fact, the reception of the logistic part of Lorenzen's article takes place mostly indirectly, through the articles \citealt{schuette51,schuette52} and the book \citealt{schuette60}:\footnote{See \S\,18 and ch.~IX\@. In a letter to Bernays dated 7 March 1957, Schütte tells the following foremost reason to write his book: `After Mister Lorenzen has published a book from his point of view [\citealt{lorenzen55}], it seems necessary to me that the axiomatic direction also has its say. At any rate, I have the impression that the Americans and also the Münster school do not rightly take notice of the results of Mister Ackermann and myself.'\deutschmark\ Note that \citealt{schuette77} is not a mere translation of \citealt{schuette60}, as the author abandons the treatment of ramified type theory in this second edition; in doing so, he conceals Lorenzen's contributions to proof theory but for a spurious mention of \citealt{lorenzen51} in the bibliography. Also the survey article \citealt{schuetteschwichtenberg90} records Lorenzen's contribution to logic in an elusive way.}\deutschtext{`Nachdem Herr Lorenzen ein Buch von seinem Standpunkt aus herausgebracht hat, erscheint es mir nötig, daß auch die axiomatische Richtung zu Wort kommt. Ich habe jedenfalls den Eindruck, daß die Amerikaner und auch die Münstersche Schule die Ergebnisse von Herrn Ackermann und mir nicht recht beachten' (Hs.~975:4234).} see e.g.\ \citealt[Appendix]{mendelson64}, \citealt{tait68}, \citealt[ch.~6]{girard87}, and \citealt[§~2.1]{girard00}.

In his book, \citet{schuette60} introduces the \textomega-rule as `rule UJ*', where `UJ' stands for `infinite induction', with a description of its meaning in terms of content by a reference to constructiveness: `For the application of rule UJ* requires a metalogical investigation. It will be requested that infinitely many formulas~\(F(z_1,\dots,z_n)\) have been proved to be derivable on the basis of general considerations before it is allowed to infer the derivability of the formula~\(F(a_1,\dots,a_n)\).' Here the~\(z_i\)'s are numerals while the \(a_i\)'s are free variables. He writes further: `We request of a metalogical proof needed for the application of the infinite induction that it is led as all metalogical investigations in a constructive way. That is: the metalogical proof is to consist in the specification of a general procedure after which the requested derivations resp.\ derivation parts may be immediately exhibited' (compare also \citealp{schuette51}, p.~369). He coins the expression `semi-formal system' for a formal system extended with the \textomega-rule. In contrast, in the second edition of his book, \citet{schuette77} keeps silent about the meaning of the \textomega-rule and states only this: `We call the system DA* \emph{semi-formal} since, as opposed to formal systems, it has basic inferences (S2.0*) with infinitely many premises' (p.~174). This silence is in stark opposition to the introduction which insists in the same terms as the first edition on the `constructive character' and the `constructive standpoint' as the framework of metamathematics. Thus the way of dealing with infinitely many premisses is considered as a private business of the proof theorist: he should not need to express on the record the meaning of infinitary proof objects and might e.g.\ resort to set theory for these. The metalogical investigation is eluded.

\citealt{schuette62}, written for a general audience, compares three methods for proving the consistency of arithmetic: Gentzen's use of transfinite induction; Lorenzen's `semi-formal' use of the \textomega-rule; Gödel's use of computable functionals of finite type. The last is the `most direct' as it possesses `a character of immediate evidence' whereas `ordinal inductions appear as admissible only after a corresponding foundation'; `semi-formal systems permit analyses of the logico-mathematical deduction that suggest themselves and are particularly transparent'; `transfinite ordinals [\dots{}] 
give us the possibility to characterise the different induction principles used metamathematically with respect to logical strength by equivalent ordinal inductions' (pp.~106--107).

The detour via Schütte's reception may have contributed to proof theory continuing to focus on measuring logical strength by ordinal numbers, whereas the fact that Lorenzen does not resort to ordinals in his proof of consistency should be considered as a feature of his approach.

\citet{tait68} provides a very clear presentation of Schütte's approach with a mention of Lorenzen; see also \citealt[§~3.2,]{fefermansieg81} for an account of \citealt{tait68}.

\subsection{The generalised inductive definition of accessibility.}
\label{sec:ackermann}

Wilhelm Ackermann hears about Lorenzen's proof of consistency in 1946 through Bernays (see his letter to Lorenzen dated 11 November 1946, PL~1-1-125). As he writes in a subsequent letter dated 21 May 1947, Ackermann is working at the time at setting up `mathematics out of a type-free logical system that is demonstrably free from contradiction'\deutsch{`aus einem nachweislich widerspruchsfreien typenfreien logischen Axiomensystem die Mathematik aufzubauen' (PL~1-1-117).} and for this he needs a `considerably higher ordinal number' than in the `transfinite inductions up to the first \textepsilon-number' that `Gentzen and [he] use in [their] proofs of freedom from contradiction of arithmetic'. He is therefore interested in `constructively recordable numbers of the second number class' following \citet{church38}.\deutsch{`So benutzen Gentzen und ich bei unseren Widerspruchsfreiheitsbeweisen für die Arithmetik transfinite Induktionen bis zur ersten \textepsilon-Zahl. Bei den Untersuchungen, an denen ich augenblicklich arbeite, gehe ich bis zu einer wesentlich höheren Ordinalzahl. Unter einer konstruktiv erfassbaren Zahl der II.~Zahlenklasse verstehe ich dabei im Anschluss an A. Church [\dots].'}

In their correspondence, they are also interested in a set-up of ordinal numbers by  Lorenzen: see Ackermann's letters of 21 May and of 9 June 1947. In the latter he writes: `As far as I understand the train of thought of your set-up from your hints, your set-up of the ordinal numbers of the 2.~number class appears to me as constructivistically usable'.\deutsch{`Soweit ich den Gedankengang Ihres Aufbaus nach Ihren Andeutungen verstehe, erscheint mir Ihr Aufbau für die Ordinalzahlen der 2.~Zahlenklasse als konstruktivistisch brauchbar' (PL~1-1-115).}

Recall that \citet[§~15.4]{gentzen36} proves the finiteness of his reduction procedure by a transfinite induction described as the generalised inductive definition of the `accessibility [Erreichbarkeit]' of the first \textepsilon-number. He discusses the constructiveness of this concept in~§~16.1\,1, which may be seen as the very heart of his proof of the consistency of elementary number theory: `[Accessibility] acquires a \emph{sense} merely by being predicated of a definite [ordinal] number for which \emph{its validity} is simultaneously \emph{proved.}'

\citet{ackermann51} sets up the segment of ordinal numbers that he needs and gives a remarkably precise account of the constructiveness of accessibility, relying in particular on Lorenzen's \textomega-rule. He provides in fact a description for the generalised inductive definition of the concept of accessibility, and we guess that it has benefited from his exchanges with Lorenzen. The `o-numbers' below are given as a certain recursively defined totally ordered system of symbols.

\begin{quotation}
  \noindent Now we have to show that the o-numbers are really ordinal numbers, or, in other words, that we are allowed to apply \emph{deductions by transfinite induction}.

  In order to lay out the text of the following considerations with less difficulty, it is advisable to introduce the following symbols, but with which we only express concisely contentful states of affairs. These symbols are: \(\mathfrak  A(\alpha)\), `the property~\(\mathfrak A\) applies to~\(\alpha\)'; \(\mathfrak K_x(\alpha,\mathfrak A(x))\), `the property~\(\mathfrak A\) applies to all o-numbers which are less than~\(\alpha\)'; \(\mathfrak V_x(\alpha,\mathfrak A(x))\), `with \(\mathfrak K_x(\beta,\mathfrak A(x))\) also \(\mathfrak A(\beta)\) has always to be the case, provided \(\beta\leqq\alpha\)'. It may now be that \(\mathfrak K_x(\alpha+1,\mathfrak A(x))\) is derivable from the assumptions~\(\mathfrak A(1)\) and~\(\mathfrak V_x(\alpha,\mathfrak A(x))\) by constructive deductions without assuming anything else about~\(\mathfrak A\). We would then first say that the number~\(\alpha\) is accessible through~\(\mathfrak A\). Now, if an o-number is accessible through~\(\mathfrak A\), then it is also accessible through any other property~\(\mathfrak B\). For, as we had not assumed anything about~\(\mathfrak A\) in the derivation of \(\mathfrak K_x(\alpha+1,\mathfrak A(x))\) but that \(\mathfrak A(1)\) and \(\mathfrak V_x(\alpha,\mathfrak A(x))\) has to be the case, so we need only to replace everywhere~\(\mathfrak A\) by~\(\mathfrak B\) in the deductions that lead from~\(\mathfrak A(1)\) and \(\mathfrak V_x(\alpha,\mathfrak A(x))\) to \(\mathfrak K_x(\alpha+1,\mathfrak A(x))\). We may therefore, instead of saying that \(\alpha\)~is accessible through~\(\mathfrak A\), simply say that \(\alpha\)~is accessible. It might now seem that accessibility is defined by a claim of generality over predicates. But this is not our conception. We want to conceive the accessibility of a number~\(\alpha\) as a certain intuitive fact, viz.\ precisely as the presence of a certain system of deductions that leads from~\(\mathfrak A(1)\) and \(\mathfrak V_x(\alpha,\mathfrak A(x))\) to \(\mathfrak K_x(\alpha+1,\mathfrak A(x))\). All deductions of the so-called positive logic are to belong to these deductions, further also the number-theoretic induction and the corresponding operating with the universal sign for o-numbers as it would e.g.\ fit intuitionistic number theory. We do not want to specify here these deductions in detail because the following proof shows which are needed. We remark only that the following deduction is also needed: if a claim may be derived for each concrete o-number, then also the corresponding universal claim is to be derivable.\,\up{2} For the use of the universal sign we remark that the o-numbers represent a countable set that is precisely defined and set up constructively, so that the use of a universal sign for o-numbers is legitimated in the same way as that of the universal sign for natural numbers.

\vspace*{2.6pt}\footnotesize\noindent\rule{2in}{0.4pt}\vspace*{2.6pt}

{\rule{0pt}{\footnotesep}2. This deduction appears self-evident from the point of view of content. It has been applied first by P. Lorenzen inside a logistic system. Cf.\ P. Lorenzen, The freedom from contradiction of the logic of ramified type (Journal of Symbolic Logic, to appear shortly).}
\end{quotation}
In a letter to Ackermann dated 3 March 1951 \german{\citep[see][p.~197]{ackermann83}}, Lorenzen writes:
\begin{quotation}
  \noindent Thank you very much for your work `Constructive set-up of a segment of the 2.~number class'---your construction impresses me very much, I have tried earlier a similar one but not as far-reaching. I wholly share your views on the constructiveness of your definitions and proofs.
\end{quotation}
\citet[§~4]{schuette54} takes over the argument of \citealt{ackermann51} and extends its system of o-numbers into so-called `Klammersymbole' that generalise also the system of ordinal fixed points of \citet{veblen08}. This argument may also be found in
\citealt[§~11--12,]{schuette60} that presents an intermediate system inspired by the coding as integers of \citealt[§~5.3.c]{hilbertbernays39}. In contrast, it is absent from \citealt{schuette77}, as is any description of the nature of a constructive argument; at the beginning of §~24 introducing higher ordinals on p.~221, the reading, constructive or axiomatic, is up to the reader: `we use the notions \emph{map} (\emph{function}) and \emph{set} in a naive way. But these may also be regarded as being determined axiomatically (in the context of a general axiom system for set theory).'

See also the account of accessibility by \citet[note~c on p.~272]{goedel90}.

\subsection{Are infinitary inductive definitions predicative?}

\label{sec:citealtschuette62}

A few years later, \citet[p.~110]{schuette62} dissociates himself from the analysis of \citealt{ackermann51} and
qualifies the generalised inductive definition of accessibility as impredicative:
\begin{quotation}
  \noindent Hereby one proceeds in an impredicative way by including the concept of accessibility itself, defined with reference to the totality of all properties of certain ordinals, into these properties.
\end{quotation}
\citet[p.~280]{schuette65} defines `an ordering relation~\(\prec\) of equivalence classes of natural numbers representing a sufficiently large segment of the second number class in a constructive way'. On p.~286, he writes: `The relation~\(\prec\) can be proved by \emph{impredicative methods} to be a \emph{well-ordering} using a proof similar to that for the related \(\prec\)-relation in §~12 of \citealt{schuette60}.\/'

In particular, Lorenzen's theorem induction would be impredicative if one followed Schütte's narrow acceptation.
  
\subsection{Further accounts.}

\citealt{lorenzen55} expands on the rôle of lattices in logic (§~7) and mathematics (part~III). \citealt{lorenzen58, lorenzen87} provide a proof of Gentzen's subformula theorem by the method of his article. \citealt{lorenzen62} (\S~7) returns to the subject of proofs of consistency.

Oskar \citet{becker54} refers to Lorenzen's article in the last pages of his book \emph{Grundlagen der Mathematik in geschichtlicher Entwicklung}.

Evert W. \citet[p.~253]{beth59} gives a short and precise account of Lorenzen's article.

Maurice \citet[p.~62--63]{meigne59} does so as well, but does not seem to understand the content of the \textomega-rule.

Haskell B. \citet[ch.~4, theorem~B9]{curry63} follows Lorenzen in characterising a distributive lattice as a lattice satisfying cut.

Manfred E.~\citet[p.~12--13]{gen1969} comments on the relationship of Lorenzen's article with Gentzen's work.

Its philosophical significance is addressed by Matthias \citet{wille13,wille16}.

\citealt{coquand21} describes Lorenzen's standpoint in a broader context than given here.

\section{Our reception: impredicative quantification and inductive de\-fi\-ni\-tions.}
\label{sec:impr-quant-induct}

\citet[p.~57]{whiteheadrussell25} see the axiom of reducibility as a generalisation of Leibniz's identity of indiscernibles. For instance, the
 formula $\forall_X~X(3)\rightarrow X(x)$, which seems impredicative since it contains a quantification over all predicates, is actually equivalent to
 the predicative formula $3 = x$. Using in a crucial way ideas from \citealt{lorenzen20},
 we extend this remark to a predicative interpretation of some
 formulae involving a seemingly impredicative universal quantification over all predicates.

 \subsection{Free Lorenzen algebras and countable choice.}

 Let us state the definition of a countably complete Lorenzen algebra.
\begin{definition}\label{lorenzen-algebra}
 A \emph{Lorenzen algebra} is a meet-semilattice with a \emph{negation}, i.e.\ a least
 element~$0$
 and an operation~$\psc a$ such that $b \leqslant \psc a$ if, and only if,
 $b\land a = 0$. (This semilattice has then automatically a greatest element $1 = \psc{0}$.) Such an algebra is \emph{countably complete} if it is endowed with an infinitary meet operation over any sequence of elements indexed by~\(\nats\).
\end{definition}

 The main result of \citealt{lorenzen20} is to essentially build the free
 countably complete Lorenzen algebra $K$ over a given preordered set $P$ and to show that the canonical
 map $P\rightarrow K$ is an embedding. This is a purely semilattice-theoretic reformulation of
 Gentzen's 1936 consistency proof of arithmetic. We write `essentially' since the actual
 statement is a little more complex if one wants to avoid the use of the axiom of countable
 choice. We provide such a statement below; note that a possible constructive
 way to avoid the axiom of countable choice is provided by the setting of \cite{hottbook}.

 \subsection{An impredicative formal system: syntax.}

 We consider the following language.
 The terms are of the form $S^k(0)$ and $S^k(x)$. A closed term $t$ represents such a natural
 number, that we will also write $t$. The atomic formulae are of the form $X(t)$
 or $P(t_1,\dots,t_n)$, where $P$ represents some $n$-ary boolean-valued function.
 
 We have formulae built from $\top,~\neg \psi,~\varphi\land \psi$,
and $\forall_x \psi$ and $\forall_X \psi$. 

A formula is arithmetical if it does not contain any quantification over predicates.

A formula is \emph{strictly $\Pi^1_1$} if it is of the form $\forall_X\mu$, where $\mu$ is
arithmetical and uses only $X$ as a predicate variable.

We consider the fragment of the language in which we form only strictly $\Pi_1^1$
universal quantifications $\forall_X \mu$.

We also have terms for predicates $T,U,\dots$. They are of the form
$T = \lambda_x\mu$. We define the substitution $\psi(T/X)$ for a closed predicate $T$
by induction on $\psi$:
\[\begin{aligned}
    (X(t))(T/X) &= \mu(t/x)&
     (\neg \mu)(T/X) &= \neg {(\mu(T/X))}\\
     (Y(t))(T/X) &= Y(t)&
     (\mu_0\land \mu_1)(T/X) &= \mu_0(T/X)\land \mu_1(T/X)\\
     P(t_1,\dots,t_n)(T/X) &= P(t_1,\dots,t_n)&
     (\forall_x\mu)(T/X) &= \forall_x (\mu(T/X))
\end{aligned}\]

 \subsection{A semilattice defined in a predicative metatheory.}

 We use \citealt{lorenzen20} to build, in a predicative metatheory, a special semilattice~\(L\).

 We let $n,m,\dots$ range over natural numbers.

 We start from an infinite set of symbols $X,Y,\dots$
 representing predicate variables on numbers.

 Then we  build inductively the symbolic objects in~$L$ by the rule
 $$
 a, b, \dots~::=~X(n)~|~0~|~\neg {a}~|~a\land b~|~\bigland_n a_n
 $$
 with~\(|\) delimiting the alternatives of the rule.
  
One way to interpret \citealt{lorenzen20} is that he defines a \emph{preorder} relation
$a\leqslant b$ on the set $L$ 
such that $a\leqslant\neg b$ if, and only if, $a\land b \leqslant 0$;
$0\leqslant a$ for all $a$; $b\leqslant \bigland_n a_n$ if, and only if, $b\leqslant a_n$ for all $n$;
$c\leqslant a\land b$ if, and only if, $c\leqslant a$ and $c\leqslant b$.

If we wanted to build the free countably complete Lorenzen algebra over the atoms $X(n)$, we would
have to quotient by the equivalence relation defined by $a\leqslant b$ and $b\leqslant a$,
and countable choice would be needed to show that we get a countably complete Lorenzen algebra.

If $\Omega$ is a countably complete Lorenzen algebra and $\rho$ assigns to any
predicate symbol~$X$ occurring in $a$ a function $\rho(X)\colon\nats\rightarrow \Omega$, we can
compute $a\rho$ in $\Omega$ by induction on~$a$ by the following rules:
\[X(n)\rho = \rho(X)(n)\quad0\rho = 0\quad(\neg {a})\rho = \pr{a\rho}\quad(a\land b)\rho = a\rho\land b\rho\quad(\bigland_n a_n)\rho = \bigland_n (a_n\rho)\]
We have $a\rho\leqslant b\rho$ in $\Omega$ whenever $a\leqslant b$ in~\(L\).

\medskip

 We let $L_{\mathrm{fin}}$ be the subset of symbolic objects in~$L$ which depend on
\emph{finitely many} predicate symbols.

An \emph{ideal} $A$ is a subset of $L_{\mathrm{fin}}$ containing $0$
and such that $b\in A$ whenever $a\in A$ and $b\leqslant a$. In particular any element
$a$ in $L_{\mathrm{fin}}$ defines the principal ideal $\down a$ of all elements in $L_{\mathrm{fin}}$ such
that $b\leqslant a$.

We let $\Omega$ be the set of all ideals.
$\Omega$ has a structure of countably complete Lorenzen algebra, with intersection as
 meet and with least element $\{0\}$.
 The negation operation defines $\pr{A}$ as the ideal
 of elements $b$ in $L_{\mathrm{fin}}$ such that $b\land a\leqslant 0$ whenever $a$ is in $A$.

 In particular, if $\rho$ assigns a function $\nats\to\Omega$ to any predicate symbol
 $X$, we can compute $a\rho$ in $\Omega$ for $a$ in $L_{\mathrm{fin}}$.

 By induction on $a$, we can show the following result.

 \begin{lemma}\label{Yoneda}
   If $a$ is in~$L_{\mathrm{fin}}$ and we have $\rho(X)(n) = \down (X(n))$ for the predicate symbols~$X$ occurring in~$a$, then 
   $a\rho = \down a$.
 \end{lemma}

 Let $c(X)$ be a symbolic object in $L_{\mathrm{fin}}$
 in which at most $X$ occurs as predicate symbol. We shall write \(c(Y)\) for the object arising by substitution of~\(X\) by~\(Y\) in~\(c(X)\), and ${c(X)(X=f)}$ for the element
 \(c(X)\rho\) computed in $\Omega$ with \(\rho\)
 assigning the function $f\colon\nats\rightarrow \Omega$ to~\(X\).
The key result which provides a predicative analysis of (strict) impredicative
 comprehension is the following.

 \begin{theorem}\label{main}
   The
   family $c(X)(X=f)$, where $f$ ranges over the functions $\nats\rightarrow\Omega$, has a g.l.b.\
   $\bigland_f c(X)(X=f)$ in the semilattice $\Omega$.
 \end{theorem}
 
 \begin{proof}
   Define $A$ to be the set of all elements $a$ such that
   $a\leqslant c(X)$ if $X$ does not occur in $a$.
   Note that if neither~$X$ nor~$Y$ occur in $a$, then $a\leqslant c(X)$ is equivalent
   to $a\leqslant c(Y)$. It follows from this remark that $a$
   is in $A$ if, and only if,
   $a\leqslant c(X)$ for \emph{all} $X$ not occurring in $a$ if, and only if, $a\leqslant c(X)$
   for \emph{some} $X$ not occurring in $a$. If $a$ is in $A$
   and $b\leqslant a$, then $a\leqslant c(X)$ for some $X$ not occurring in $a$; if we consider $Y$ occurring \emph{neither} in~$a$ \emph{nor} in~$b$,
   then $a\leqslant c(Y)$ and hence $b\leqslant c(Y)$,
   so that $b$ is in $A$.
   The set $A$ is thus an ideal, i.e.\ an element of $\Omega$.

   We claim that $a$ is in $A$ if, and only if, it belongs to all $c(X)(X=f)$, which
   will show that $A$ is the g.l.b.\ of the family $c(X)(X=f)$.

   If $a$ is in $A$ then $a\leqslant c(X)$ for $X$ not occurring in $a$.
   Define $\rho(Z)(n) = \down (Z(n))$ for $Z$ occurring in $a$ and $\rho(X) = f$.
   We then have $a\rho\leqslant  c(X)\rho$.
   But $a\rho = \down a$ by lemma~\ref{Yoneda}
   and $c(X)\rho = c(X)(X = f)$. Hence $\down a \leqslant c(X)(X=f)$
   and $a$ is in $c(X)(X=f)$.

   If conversely $a$ contains no occurrence of~$X$ and is in all $c(X)(X = f)$, where $f$ ranges over~$\nats\rightarrow \Omega$, then in particular it is in $c(X)(X=g)$ for
   $g(n) = \down (X(n))$. In this case $c(X)(X=g)$ is $\down c(X)$ by lemma~\ref{Yoneda} and so $a\leqslant c(X)$,  i.e.\ $a$ is in $A$.
 \end{proof}
 
 \subsection{Interpretation of strict \texorpdfstring{$\Pi^1_1$}{Pi 1 1}-comprehension: semantics.}

Any arithmetical formula $\nu$ defines an element $[\nu]$ in $L_{\mathrm{fin}}$ by the rules
\[
  \begin{gathered}
\relax[X(t)] = X(t)\qquad
[P(t_1,\dots,t_n)] = \delta_P(t_1,\dots,t_n)\\
[\neg \mu] = \neg {[\mu]}\qquad[\mu_0\land \mu_1] = [\mu_0]\land [\mu_1]\qquad
[\forall_x\mu] = \bigland_n [\mu(n/x)]
\end{gathered}
\]
where \(\delta_P\) is the $n$-ary \(0,1\)-valued function representing~\(P\). 

If $\rho$ assigns a function $\rho(X)\colon\nats\rightarrow\Omega$
for $X$ free in $\psi$, we can define
the semantics $\sem{\psi}\rho$ as an element of $\Omega$.
We define it first for an arithmetical formula by the clauses
\[
  \begin{gathered}
\sem{X(t)}\rho = \rho(X)(t)\qquad
\sem{P(t_1,\dots,t_n)} = \delta_P(t_1,\dots,t_n)\\
\sem{\neg \mu}\rho = \neg{(\sem{\mu}\rho)}\qquad
\sem{\mu_0\land \mu_1}\rho = \sem{\mu_0}\rho\cap \sem{\mu_1}\rho\qquad
\sem{\forall_x\mu}\rho = \bigcap\nolimits_n \sem{\mu(n/x)}\rho
\end{gathered}
\]
We can now define $\sem{T}\rho$ for $T = \lambda_x\mu$ by
$\sem{\lambda_x\mu}\rho (n) = \sem{\mu(n/x)}\rho$.

 If $\mu$ is an arithmetical formula with at most one free variable $X$, then
 $\sem{\mu}(X=f)$ is equal to $[\mu](X=f)$.
 By theorem~\ref{main}, the family $\sem{\mu}(X=f)$ has a g.l.b.\ and
 we can define
\[\sem{\forall_X \mu}\rho = \bigland_{f} \sem{\mu}(\rho,X= f)\text.\]
 In particular, we get
\[\sem{\forall_X\mu}\rho \leqslant \sem{\mu}(\rho, X = \sem{T}\rho) = \sem{\mu(T/X)}\rho\]
 and our semantics, which we have built in a predicative metatheory, justifies
 comprehension for strictly $\Pi_1^1$ formulae.


\def\Dbar{\leavevmode\lower.6ex\hbox to 0pt{\hskip-.23ex \accent"16\hss}D}

\end{document}